\documentclass[3p,review]{elsarticle}

\usepackage{amssymb}  
\usepackage{amsthm}   
\usepackage{lineno}   
\usepackage{amsmath}  
\usepackage{cases}    
\usepackage{tabularx} 
\usepackage[hidelinks]{hyperref} 
\hypersetup{pdfstartview={XYZ null null 1.00}, pdfpagemode=UseNone}
\usepackage[figuresright]{rotating}  
\usepackage{color}

\begin{document}

\begin{frontmatter}


\title{A Simple, Efficient, and High-order Accurate Curved Sliding-mesh Interface Approach to Spectral Difference Method on Coupled Rotating and Stationary Domains}

\author[GWU]{Bin Zhang}

\author[GWU]{Chunlei Liang \corref{cor1}}
\ead{bzh@gwu.edu,chliang@gwu.edu}

\cortext[cor1]{Corresponding author. Tel: +1 (202)-994-7073}

\address[GWU]{Department of Mechanical and Aerospace Engineering, \\ The George Washington University, Washington, DC 20052}

\begin{abstract}
This paper presents a simple, efficient, and high-order accurate sliding-mesh interface approach to the spectral difference (SD) method. We demonstrate the approach by solving the two-dimensional compressible Navier-Stokes equations on quadrilateral grids. This approach is an extension of the straight mortar method originally designed for stationary domains \cite{kopriva-96-jcp2,kopriva98}, it employs curved dynamic mortars on sliding-mesh interfaces to couple rotating and stationary domains. On the nonconforming sliding-mesh interfaces, the related variables are first projected from cell faces to mortars to compute common fluxes, and then the common fluxes are projected back from the mortars to the cell faces to ensure conservation. To verify the spatial order of accuracy of the sliding-mesh spectral difference (SSD) method, both inviscid and viscous flow cases are tested. It is shown that the SSD method preserves the high-order accuracy of the SD method. Meanwhile, the SSD method is found to be very efficient in terms of computational cost. This novel sliding-mesh interface method is very suitable for parallel processing with domain decomposition. It can be applied to a wide range of problems, such as the aerodynamics of rotorcraft, wind turbines, and oscillating wing wind power generators, etc.
\end{abstract}

\begin{keyword}
spectral difference method \sep rotating mesh \sep sliding mesh \sep high-order methods \sep unstructured grid
\end{keyword}

\end{frontmatter}



\section{Introduction}
\label{sec:introduction}

High-order (third and above) numerical methods are becoming more and more popular in recent years due to their capability of producing more accurate solutions on relatively coarse grid \cite{ZJWang-14}. The spectral difference (SD) method is one discontinuous high-order method for solving the conservation laws on unstructured grids \cite{liu-jcp-2006,wang-jsc-2007,sun-2007-ccp,liang-caf-2009}. This method is an extension of the staggered multi-domain high-order method originally designed by \citet{kopriva-96-jcp1}. It was shown that the SD method also has strong connection with the Flux Reconstruction/Correction Procedure via Reconstruction (FR/CPR) methods \cite{Huynh-14}, and it shares similarity with the quadrature-free discontinuous Galerkin method \cite{May-11}. The stability of a particular choice of flux points for the SD method was proved by \citet{jameson-2010-jsc} for the one-dimensional linear wave equation. Although the proof has not been generalized to higher-dimensional tensor-product elements, we have not observed numerical instability from several successful turbulent flow simulations \cite{Liang-09,Castonguay-01,Abrar-10}.

We have seen more and more applications of the SD method to realistic flow simulations, for example, for large eddy simulations on fixed grids \cite{Liang-09,Castonguay-01,Abrar-10,Parsani-10,Parsani-11,Lodato-14}. The SD method is also particularly promising for simulating vortex-dominated flows on moving and deforming grids \cite{Ou-10,Yu-2011-caf}. \citet{liang-2011-caf} extended the SD method for simulating two-dimensional unsteady flows around a plunging or pitching airfoil. \citet{Dejong-14} studied three-dimensional vortex induced vibrations using the SD method. 

However, when the mesh undergoes very large rotation motion, such as for flows around rotating propellers or passing a flapping wing with very large pitching angles, remeshing \cite{Tezduyar-1992a,Tezduyar-1992b} is required. Our goal is to involve the minimum number of remeshing and simultaneously preserve the high-order accuracy of the SD method. This motivates us to develop a new approach to the SD method for coupled rotating and stationary domains with sliding-mesh interfaces. In our approach, both inviscid and viscous fluxes on the sliding-mesh interfaces are constructed using a newly developed curved dynamic mortar method. The mortar method on fixed grids was originally proposed for incompressible flows by \citet{mavriplis1989}. \citet{kopriva-96-jcp2,kopriva98} proved the conservation property of the mortar method for the compressible flow equations and applied it to the compressible Euler and Navier-Stokes equations on stationary domains using structured grids. In this paper, we show that our sliding-mesh approach is as simple as those designed for low-order numerical methods \cite{Yeoh-04,McNaughton-14} while preserving the high-order accuracy of the SD method. This simple but novel sliding-mesh spectral difference (SSD) method can have a wide range of applications, such as rotorcraft aerodynamics, wind turbine wake dynamics, and oscillating wing wind power generators.

The paper is organized as follows: Section \ref{sec:equations} gives the two-dimensional compressible Navier-Stokes equations on stationary and rotating domains. Section \ref{sec:methods} reviews the SD method and presents the SSD method in detail. Verification studies and applications are reported in Section \ref{sec:tests}. Section \ref{sec:conclusion} concludes the paper.

\section{The governing equations}
\label{sec:equations}

\subsection{The compressible Navier-Stokes equations on stationary domain}

We consider the two-dimensional unsteady compressible Navier-Stokes equations in conservative form,
\begin{equation}
    \frac{\partial \mathbf{Q}} {\partial t} +  \frac{\partial \mathbf{F}} {\partial x} + \frac{\partial \mathbf{G}} {\partial y} = 0,
    \label{eq:physical}
\end{equation}
where $\mathbf{Q}$ is the vector of conservative variables, $\mathbf{F}$ and $\mathbf{G}$ are the $x$ and the $y$ flux vectors. These terms have the following expressions,
\begin{align}
    \mathbf{Q} &= [\rho ~ \rho u ~ \rho v ~ E]^T, \label{eq:Q} \\
    \mathbf{F} &= \mathbf{F}_{inv}(Q) + \mathbf{F}_{vis}(Q,\nabla Q), \label{eq:F}\\
    \mathbf{G} &= \mathbf{G}_{inv}(Q) + \mathbf{G}_{vis}(Q,\nabla Q), \label{eq:G}
\end{align}
where $\rho$ is the fluid density, $u$ and $v$ are the $x$ and the $y$ velocities, $E$ is the total energy per volume defined as $E = p/(\gamma-1) + \frac{1}{2}\rho(u^2+v^2)$, $p$ is the pressure, $\gamma$ is the ratio of specific heats and is set to 1.4 (i.e., the typical value for the air in
standard conditions).

As shown in Equations (\ref{eq:F}) and (\ref{eq:G}), the fluxes have been divided into inviscid and viscous parts. The inviscid fluxes are only functions of conservative variables, which are
\begin{equation}
    \mathbf{F}_{inv} = \left[
    \begin{array}{c}
         \rho u       \\
         \rho u^2 + p \\
         \rho uv      \\
         (E+p)u
    \end {array}
    \right], \
    \mathbf{G}_{inv} = \left[
    \begin{array}{c}
         \rho v       \\
         \rho uv      \\
         \rho v^2 + p \\
         (E+p)v
    \end {array}
    \right].
\end{equation}
The viscous fluxes are functions of the conservative variables as well as their gradients. They have the following expressions,
\begin{equation}
    \mathbf{F}_{vis} = -\left[
    \begin{array}{c}
         0         \\
         \tau_{xx} \\
         \tau_{yx} \\
         u\tau_{xx}+v\tau_{yx}+\mathrm{k}T_x
    \end {array}
    \right], \
    \mathbf{G}_{vis} = -\left[
    \begin{array}{c}
         0         \\
         \tau_{xy} \\
         \tau_{yy} \\
         u\tau_{xy}+v\tau_{yy}+\mathrm{k}T_y
    \end {array}
    \right],
\end{equation}
where $\tau_{ij}$ is the shear stress tensor and is related to the velocity gradients as $\tau_{ij} = \mu (u_{i,j}+u_{j,i}) + \lambda\delta_{ij}u_{k,k}$, $\mu$ is the dynamic viscosity, $\lambda=-2/3\mu$ based on Stokes' hypothesis, $\delta_{ij}$ is the Kronecker delta, $\mathrm{k}$ is the thermal conductivity, $T$ is the temperature which is related to density and pressure through the ideal gas law $p=\rho R T$, where $R$ is the gas constant.


\subsection{The compressible Navier-Stokes equations on rotating domain}

On the rotating domains, we implement a simplified equation which is equivalent to the Arbitrary Lagrange-Eulerian (ALE) \cite{Hirt-74} form of Equation (\ref{eq:physical}). Due to grid motion, the inviscid fluxes are modified to take the following forms,
\begin{equation}
    \mathbf{F}_{inv} = \left[
    \begin{array}{c}
        \rho u       \\
        \rho u^2 + p \\
        \rho uv      \\
        (E+p)u
    \end {array}
    \right] - u_g\left[
    \begin{array}{c}
        \rho    \\
        \rho u  \\
        \rho v  \\
        E
    \end {array}
    \right], \
    \mathbf{G}_{inv} = \left[
    \begin{array}{c}
        \rho v       \\
        \rho uv      \\
        \rho v^2 + p \\
        (E+p)v
    \end {array}
    \right] - v_g\left[
    \begin{array}{c}
        \rho    \\
        \rho u  \\
        \rho v  \\
        E
    \end {array}
    \right],
\end{equation}
where $u_g$ and $v_g$ are the $x$ and the $y$ grid velocities, respectively. The viscous fluxes and all other variables are uninfluenced and take the same expressions as those in the previous section.

For a domain rotating at angular velocity $\boldsymbol{\omega}$, the grid velocities are $(u_g,v_g) = \boldsymbol{\omega} \times \mathbf{r}$, where $\mathbf{r}$ is the position vector with respect to rotating center. For all test cases in the present study, $\boldsymbol{\omega}$ is known as a priori, thus grid velocities and coordinates are updated analytically on the rotating domains.

\subsection{The transformed equations}
As will be discussed in the next section, we map each quadrilateral grid cell in the physical domain to a standard square element in a computational domain. This mapping facilitates the construction of solution and flux polynomials. As a result, we only need to solve a set of transformed equations within each standard element. Let us assume that the physical coordinates $(x,y)$ are mapped to the computational ones $(\xi , \eta)$ through a transformation: $x=x(\xi,\eta)$, $y=y(\xi,\eta)$. It can be shown that Equation (\ref{eq:physical}) will take the following conservative form after coordinates transformation,
\begin{equation}
    \frac{\partial \widetilde{\mathbf{Q}}} {\partial t} +  \frac{\partial \widetilde{\mathbf{F}}} {\partial \xi} +
    \frac{\partial \widetilde{\mathbf{G}}} {\partial \eta} = 0,
    \label{eq:computational}
\end{equation}
where $\widetilde{\mathbf{Q}}=|\mathcal{J}|\mathbf{Q}$, and the transformed fluxes $\widetilde{\mathbf{F}}$, $\widetilde{\mathbf{G}}$ are related to the physical ones as
\begin{equation}
    \left(
    \begin{array}{c}
        \widetilde{\mathbf{F}} \\
        \widetilde{\mathbf{G}}
    \end{array}
    \right)
    = |\mathcal{J}|\mathcal{J}^{-1}
    \left(
    \begin{array}{c}
        \mathbf{F} \\
        \mathbf{G}
    \end{array}
    \right),
    \label{eq:transform_relation}
\end{equation}
where $|\mathcal{J}|$ is determinant of the Jacobian matrix, and $\mathcal{J}^{-1}$ is the inverse Jacobian matrix,
\begin{align}
    |\mathcal{J}| &= \left| \frac{\partial(x,y)}{\partial(\xi,\eta)} \right| = \left|
    \begin{array}{cc}
        x_{\xi} & x_{\eta}  \\
        y_{\xi} & y_{\eta}
    \end{array}    
    \right| = x_{\xi} y_{\eta} - x_{\eta} y_{\xi}, \\
    \mathcal{J}^{-1} &= \frac{\partial(\xi,\eta)}{\partial(x,y)} = \left[
    \begin{array}{cc}
        \xi_x  & \xi_y   \\
        \eta_x & \eta_y
    \end {array}
    \right] = \frac{1}{|\mathcal{J}|} \left[
    \begin{array}{cc}
        y_\eta & -x_\eta   \\
       -y_\xi  &  x_\xi
    \end {array}
    \right].
\end{align}

\section{Numerical methods}
\label{sec:methods}

In this section, we first give a brief review of the SD method. Subsequently, we describe a newly formulated sliding-mesh interface technique that is built on the SD formulation. For temporal discretization, an explicit strong stability preserving Runge-Kutta method \cite{ruuth-02} is used for all computations throughout this paper.

\subsection{The SD method}
For SD method on quadrilateral grids, we first transform each cell in the physical domain to a standard square element $(0\leq \xi \leq 1, 0\leq \eta \leq 1)$ in the computational domain. The transformation can be done through iso-parametric mapping. As was reported by \citet{liang-jcp,liang-2009-cas}, using linear cell defined by four nodes is not sufficient for problems involving curved boundaries. High-order cubic cell with twelve nodes are used along the curved boundaries to ensure stability and accuracy in the present study.

\begin{figure}[!ht]
    \centering
    \includegraphics[scale=1.0]{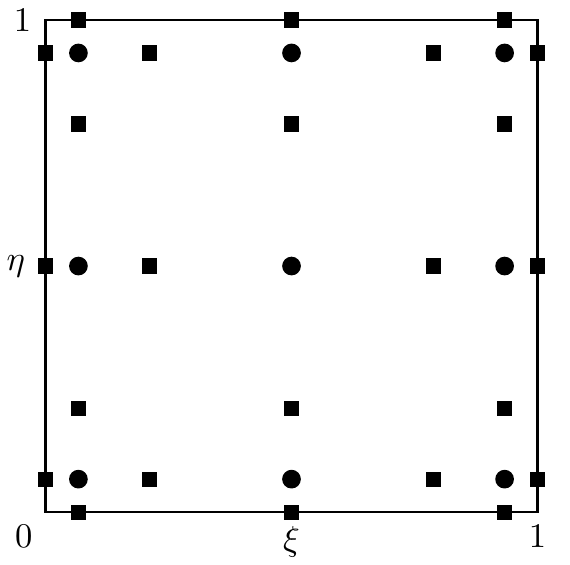}
    \caption{Schematic of the distribution of solution points (circles) and flux points (squares) for a third-order SD scheme.}
    \label{fig:spfp}
\end{figure}

After the mapping, solution points (SPs) and flux points (FPs) are defined on each standard element as shown in Figure \ref{fig:spfp} for a third-order scheme. For an $N$-$th$ order SD method, $N$ SPs are required along each coordinate direction to construct degree $(N-1)$ solution polynomials, and $(N+1)$ FPs are employed in each direction to construct degree $N$ flux polynomials. In the current implementation, the SPs: $X_s$, where $s=1,2,...,N$, are chosen as $N$ Chebyshev-Gauss points. The FPs: $X_{s+1/2}$, where $s=0,1,2,...N$, are chosen as $(N-1)$ Legendre-Gauss points plus two end points to align in a staggered fashion with the SPs.

To construct solution and flux polynomials, the following Lagrange bases at the SPs and FPs are used,
\begin{align}
    h_i(X)       &= \prod_{s=1,s\neq i}^{N}\left(\frac{X-X_s}{X_i-X_s}\right), \label{eq:hbasis} \\
    l_{i+1/2}(X) &= \prod_{s=0,s\neq i}^{N}\left(\frac{X-X_{s+1/2}}{X_{i+1/2}-X_{s+1/2}}\right).  \label{eq:lbasis}
\end{align}
The solution and the fluxes within each element are simply tensor products of the Lagrange bases,
\begin{align}
    \mathbf{Q}(\xi,\eta)             &= \sum_{j=1}^{N} \sum_{i=1}^{N} \frac{\widetilde{\mathbf{Q}}_{i,j}}{|\mathcal{J}_{i,j}|} h_i(\xi) \cdot h_j(\eta), \label{eq:Qt} \\
    \widetilde{\mathbf{F}}(\xi,\eta) &= \sum_{j=0}^{N} \sum_{i=0}^{N} \widetilde{\mathbf{F}}_{i+1/2,j} l_{i+1/2}(\xi) \cdot h_j(\eta), \label{eq:Ft}\\
    \widetilde{\mathbf{G}}(\xi,\eta) &= \sum_{j=0}^{N} \sum_{i=0}^{N} \widetilde{\mathbf{G}}_{i,j+1/2} h_i(\xi) \cdot l_{j+1/2}(\eta). \label{eq:Gt}
\end{align}

The above reconstructed solution and fluxes are only element-wise continuous, but discontinuous across cell interfaces. A Riemann solver is employed to compute the common inviscid fluxes at cell interfaces to ensure conservation. In the current implementation, the Rusanov solver \cite{rusanov-61} has been used for this purpose. The common viscous fluxes are computed from common solution and common gradients, and the detailed steps can be found in previous papers by \citet{liang-jcp,liang-2009-cas}.

\subsection{The sliding-mesh interface approach}
\begin{figure}[!ht]
    \centering
    \includegraphics[scale=1.0]{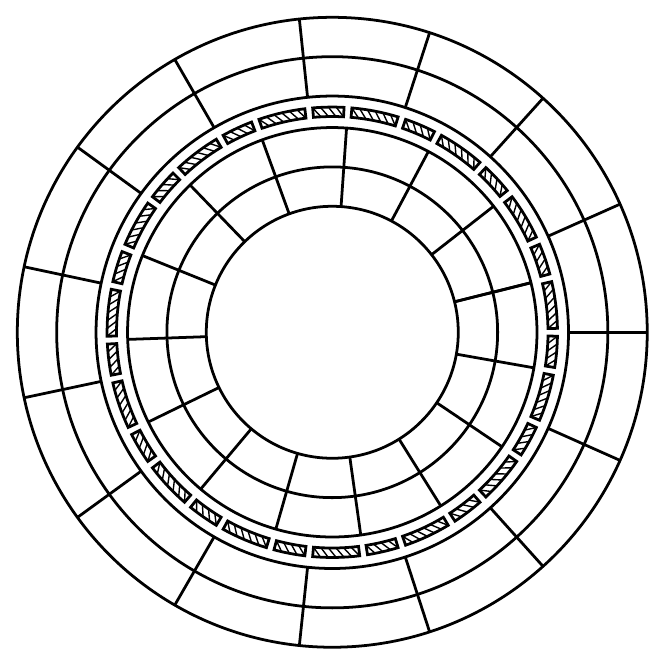} \\
    \caption{Schematic of the distribution of mortars (hatched) between a rotating and a stationary meshes.}
    \label{fig:mortar_schematic}
\end{figure}

Sliding-mesh interfaces are formed between rotating and stationary meshes. The simplest situation involves only one rotating mesh and one stationary mesh as shown in Figure \ref{fig:mortar_schematic}. The inner mesh can rotate while the outer is fixed, or vice versa. Communication between the stationary and the rotating meshes are realized through ``mortars". To make the explanation intuitive, we have scaled the inner mesh in order to place mortars in between the two coupled meshes. 

The mortars are arranged in a counterclockwise order. We refer the inner mesh as left (L) and the outer mesh as right (R) with respect to the mortars. To facilitate code implementation and reduce computational cost, cell faces on both sides of the sliding-mesh interface have been uniformly meshed. A closer look at Figure \ref{fig:mortar_schematic} reveals how mortars and cell faces on the sliding-mesh interface are connected: at each time instant, a cell face is connected to two mortars, and each mortar is associated with one left and one right cell faces. This cell face and mortar connectivity needs to be updated at every stage of the Runge-Kutta time-stepping method. As was discussed in \cite{kopriva-96-jcp2} for stationary grid, each cell face can have more than two mortars.Thus, our sliding-mesh interface method can also be extended to non-uniform meshes.

\begin{figure}[!ht]
    \centering
    \vspace{0.3cm}
    \includegraphics[scale=1.0]{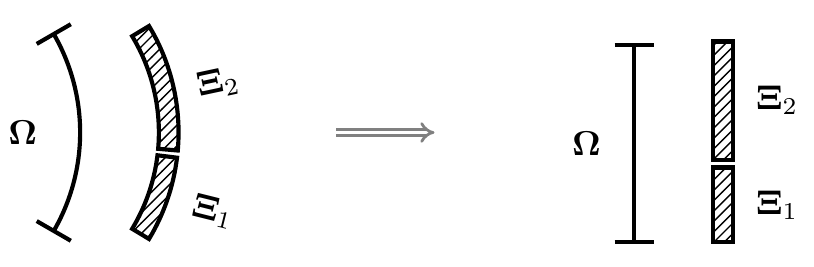} \\
    \caption{Mapping of curved cell face and mortars to straight ones: left, curved face and mortars in physical domain; right, straight face and mortars in computational domain.} 
    \label{fig:mortar_mapping}
\end{figure}

Figure \ref{fig:mortar_mapping} shows a cell face $\Omega$ and the attached two mortars $\Xi_1$ and $\Xi_2$. Each curved mortar is mapped to a straight edge $0\leq z \leq 1$ through $1D$ iso-parametric mapping. Face $\Omega$ is mapped to a straight edge $0\leq \xi \leq 1$ when the associated cell is mapped to a standard square element, thus no extra mapping is required. $\xi$ and $z$ are related by
\begin{equation}
    \xi=o(t)+s(t)z,  \label{eq:xz}
\end{equation}
where $o(t)$ is the offset of the mortar relative to the bottom node of $\Omega$ at time $t$, and $s(t)$ is the relative scaling. For the example shown in Figure \ref{fig:mortar_mapping}, we have $o_1=0$ and $s_1=L^{\Xi_1}/L^{\Omega}$ for $\Xi_1$, $o_2=L^{\Xi_1}/L^{\Omega}$ and $s_2=L^{\Xi_2}/L^{\Omega}$ for $\Xi_2$, where $L$ denotes the physical length of the face or the mortar elements.

According to Equation (\ref{eq:Qt}), the solutions on $\Omega$ can be represented as
\begin{equation}
    \mathbf{Q}^{\Omega} = \sum_{i=1}^{N}\mathbf{Q}_i^{\Omega}h_i(\xi), \label{eq:Qface}
\end{equation}
where $\mathbf{Q}_i^{\Omega}$ represents solution at the $i$-th SP on $\Omega$, and $h_i$ is the Lagrange basis defined in Equation (\ref{eq:hbasis}). If we define the same set of SPs on $0\leq z \leq 1$ for each mortar, then the solutions on each mortar element can be reconstructed similarly as
\begin{equation}
    \mathbf{Q}^{\Xi}    = \sum_{i=1}^{N}\mathbf{Q}_i^{\Xi}h_i(z), \label{eq:Qmortar}
\end{equation}
where $\mathbf{Q}_i^{\Xi}$ is the solution at the $i$-th SP on a mortar $\Xi$.

\begin{figure}[!ht]
    \vspace{0.3cm}
    \centering
    \includegraphics[scale=1.0]{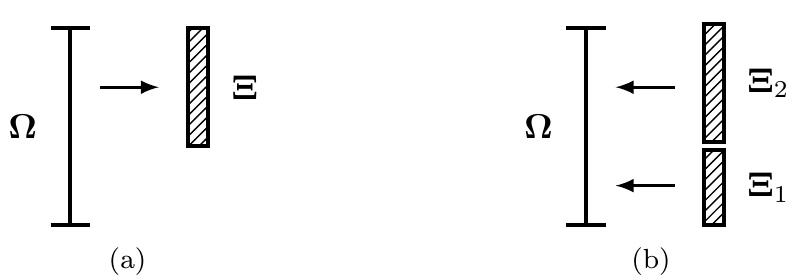} \\
    \caption{Projection between face and mortar: (a) from left face to left side of mortar, (b) from two mortars back to the associated left face.}
    \label{fig:mortar_projection}    
\end{figure}

The procedure for computing $\mathbf{Q}_i^{\Xi}$ is demonstrated in Figure \ref{fig:mortar_projection}(a). For simplicity, we only show the process on the left side of mortar $\Xi$. To get the solutions, we require that
\begin{equation}
    \int_0^1 (\mathbf{Q}^{\Xi,L}(z) - \mathbf{Q}^{\Omega}(\xi)) h_j(z)dz = 0,~~ j=1,2,...,N.
    \label{eq:prjQ1}
\end{equation}
It was shown in \cite{kopriva-96-jcp2} that the above requirement is equivalent to an unweighted $L_2$ projection. Substituting Equations (\ref{eq:xz})-(\ref{eq:Qmortar}) into the above equation and evaluating it at each SP on $\Xi$ will give a system of linear equations. The solution of this system when written in matrix form is
\begin{equation}
    \mathbf{Q}^{\Xi,L} = \mathbf{P}^{\Omega\rightarrow\Xi}\mathbf{Q}^{\Omega} = \mathbf{M}^{-1} \mathbf{S}^{\Omega\rightarrow\Xi} \mathbf{Q}^{\Omega},
    \label{eq:prjQ2}
\end{equation}
where $\mathbf{P}^{\Omega\rightarrow\Xi}$ is the projection matrix from $\Omega$ to $\Xi$, and the elements of the matrices $\mathbf{M}$ and $\mathbf{S}^{\Omega\rightarrow\Xi}$ are
\begin{align}
    M_{i,j} &= \int_0^1 h_i(z) h_j(z) dz, ~~ i,j=1,2,...,N, \label{eq:mortarM} \\
    S_{i,j}^{\Omega\rightarrow\Xi} &= \int_0^1 h_i(o+sz) h_j(z) dz, ~~ i,j=1,2,...,N, \label{eq:mortarS}
\end{align}
where $o$ and $s$ are the offset and the scaling of $\Xi$ with respect to $\Omega$. It is important to note that $o$ and $s$ are time-dependent for the sliding-mesh interface method.

The right solution vector $\mathbf{Q}^{\Xi,R}$ can be computed in the same way. Having both the left and the right solutions on a mortar, the Rusanov solver is employed to compute the common inviscid flux $\mathbf{F}_{inv}^{\Xi}$. This flux is then transformed to the computational flux as $\widetilde{\mathbf{F}}_{inv}^{\Xi}$ according to Equation (\ref{eq:transform_relation}).

As shown in Figure \ref{fig:mortar_projection}(b), to project the common inviscid fluxes $\widetilde{\mathbf{F}}_{inv}^{\Xi_1}$ and $\widetilde{\mathbf{F}}_{inv}^{\Xi_2}$ back to face $\Omega$, we require that,
\begin{equation}
    \int_0^{o_2} (\mathbf{\widetilde{F}}_{inv}^{\Omega}(\xi) - \widetilde{\mathbf{F}}_{inv}^{\Xi_1}(z)) h_j(\xi)d\xi + \int_{o_2}^1 (\widetilde{\mathbf{F}}_{inv}^{\Omega}(\xi) - \widetilde{\mathbf{F}}_{inv}^{\Xi_2}(z)) h_j(\xi)d\xi = 0, ~~ j=1,2,...,N,
    \label{eq:prjF1}
\end{equation}
where $\mathbf{\widetilde{F}}_{inv}^{\Omega}(\xi)$ is the inviscid flux polynomial on face $\Omega$. The solution of the above equation when written in matrix form is
\begin{equation}
    \widetilde{\mathbf{F}}_{inv}^{\Omega} = \mathbf{P}^{\Xi_1\rightarrow\Omega}\widetilde{\mathbf{F}}_{inv}^{\Xi_1} + \mathbf{P}^{\Xi_2\rightarrow\Omega}\widetilde{\mathbf{F}}_{inv}^{\Xi_2} = s_1\mathbf{M}^{-1}\mathbf{S}^{\Xi_1\rightarrow\Omega}\widetilde{\mathbf{F}}_{inv}^{\Xi_1} + s_2\mathbf{M}^{-1}\mathbf{S}^{\Xi_2\rightarrow\Omega}\widetilde{\mathbf{F}}_{inv}^{\Xi_2}, 
    \label{eq:prjF2}
\end{equation}
where the matrix $\mathbf{M}$ is identical to that of Equation (\ref{eq:prjQ2}), the matrices $\mathbf{S}^{\Xi_1\rightarrow\Omega}$ and $\mathbf{S}^{\Xi_2\rightarrow\Omega}$ are simply the transposes of $\mathbf{S}^{\Omega\rightarrow\Xi_1}$ and $\mathbf{S}^{\Omega\rightarrow\Xi_2}$, respectively.

For the computation of the common viscous fluxes, we first compute the common solution on each mortar as the average of the left and the right solutions,
\begin{equation}
\mathbf{Q}^{\Xi}=\frac{1}{2}(\mathbf{Q}^{\Xi,L} + \mathbf{Q}^{\Xi,R}).
\end{equation}
This common solution is then projected back to the cell faces in the same procedure as for the inviscid flux in Equation (\ref{eq:prjF2}). After that,  solution gradients and viscous fluxes are updated on the cell faces on both sides of the interface. The viscous fluxes $\widetilde{\mathbf{F}}_{vis}^{\Omega}$ on the cell faces are projected to the mortars in the same way as Equation (\ref{eq:prjQ2}). The common viscous flux $\widetilde{\mathbf{F}}_{vis}^{\Xi}$ on a mortar is taken as the average of the left and the right viscous fluxes,
\begin{equation}
\widetilde{\mathbf{F}}_{vis}^{\Xi}=\frac{1}{2}(\widetilde{\mathbf{F}}_{vis}^{\Xi,L} + \widetilde{\mathbf{F}}_{vis}^{\Xi,R}).
\end{equation}
The final step is to project $\widetilde{\mathbf{F}}_{vis}^{\Xi}$ back to cell faces, which is identical to the process showed in Equation (\ref{eq:prjF2}).

Since a uniform mesh is used for the cell faces on the sliding-mesh interface, the matrix $\mathbf{S}$ only needs to be computed for the first two mortars, and can be reused by other corresponding mortars. At the same time since the matrix $\mathbf{M}$ is time independent, it can be computed and stored before the actual calculation. To compute the integrals in Equations (\ref{eq:mortarM}) and (\ref{eq:mortarS}), one can use the Clenshaw-Curtis quadrature method as was used in \cite{kopriva-96-jcp2}. In this study, the integrand is casted into a general form as a product of $2(N-1)$ first degree polynomials, and we implement a recursive algorithm to compute the integrals analytically. This approach requires the least number of operations which is much more efficient than the Clenshaw-Curtis quadrature method.

\section{Numerical tests}
\label{sec:tests}

In this section we test the spatial accuracy of the SSD method on both inviscid and viscous flows, and then apply this method to study an external and an internal flows. A five-stage fourth-order Runge-Kutta method for time stepping \cite{ruuth-02} is used for all test cases. In each test case for demonstrating the orders of accuracy of spatial discretizations, the time step size is reduced successively until the final errors do not change with it. This ensures that the temporal discretization errors are negligible and the final errors can represent the spatial discretization errors.

\subsection{Euler vortex flow}
\label{sec:euler}

Euler vortex flow has been widely used to test the accuracies of inviscid flow solvers \cite{shu1997,ZJWang-14}. In this problem, an isentropic vortex is superimposed to and convected by a uniform mean flow. The Euler vortex flow in an infinite domain at time $t$ can be analytically described as
\begin{align}
u    &= U_{\infty}   \left\{\cos\theta - \frac{\epsilon y_r}{r_c}\exp\left(\frac{1-x_r^2-y_r^2}{2r_c^2}\right) \right\}, \\
v    &= U_{\infty}   \left\{\sin\theta + \frac{\epsilon x_r}{r_c}\exp\left(\frac{1-x_r^2-y_r^2}{2r_c^2}\right) \right\}, \\
\rho &= \rho_{\infty}\left\{1 - \frac{(\gamma-1) (\epsilon M_{\infty})^2}{2} \exp\left(\frac{1-x_r^2-y_r^2}{r_c^2}\right) \right\}^{\frac{1}{\gamma -1}}, \\
p    &= p_{\infty}   \left\{1 - \frac{(\gamma-1) (\epsilon M_{\infty})^2}{2} \exp\left(\frac{1-x_r^2-y_r^2}{r_c^2}\right) \right\}^{\frac{\gamma}{\gamma -1}},
\end{align}
where $U_{\infty}$, $\rho_{\infty}$, $p_{\infty}$, $M_{\infty}$ are the mean flow speed, density, pressure and Mach number, respectively. $\theta$ is the direction of the mean flow (i.e. the direction along which the vortex is convected), $\epsilon$ and $r_c$ can be interpreted as the vortex strength and size. The relative coordinates $(x_r,y_r)$ are defined as
\begin{align}
x_r &= x - x_0 - \bar{u} t, \label{eq:xr} \\
y_r &= y - y_0 - \bar{v} t, \label{eq:yr}
\end{align}
where $\bar{u}=U_{\infty}\cos\theta$, $\bar{v}=U_{\infty}\sin\theta$ are the $x$ and $y$ components of the mean velocity, $(x_0,y_0)$ is the initial position of the vortex. The analytical solution of the Euler vortex flow problem within a square domain ($0\leq x, y \leq L$) with periodic boundary conditions can be achieved by replacing the relative coordinates with the following expressions,
\begin{align}
x_r &= x_r - \lfloor \frac{x_r+x_0}{L} \rfloor \cdot L, \label{eq:xrnew} \\
y_r &= y_r - \lfloor \frac{y_r+y_0}{L} \rfloor \cdot L, \label{eq:yrnew}
\end{align}
where the floor operator $\lfloor x \rfloor$ gives the largest integer that is not greater than a real number $x$. The $x_r$ and $y_r$ on the right hand sides of Equations (\ref{eq:xrnew}) and (\ref{eq:yrnew}) are from Equations (\ref{eq:xr}) and (\ref{eq:yr}).

In this test, the uniform mean flow is chosen as $(U_{\infty}, \rho_{\infty}, p_{\infty}) = (1,1,1)$ with a Mach number of $M_{\infty}=0.3$. The flow direction is set to $\theta=\arctan(1/2)$. A vortex with parameters: $\epsilon=1$, $r_c=1$, is superimposed to the mean flow. The domain size is $0\leq x, y \leq 10$ (i.e. $L=10$), and the vortex is initially located at the domain center $(x_0,y_0) = (5,5)$. Periodic boundary conditions are applied in both $x$ and $y$ directions.

\begin{figure}[!ht]
\vspace{0.5cm}
  \centering
  \includegraphics[width=0.40 \textwidth, angle=0]{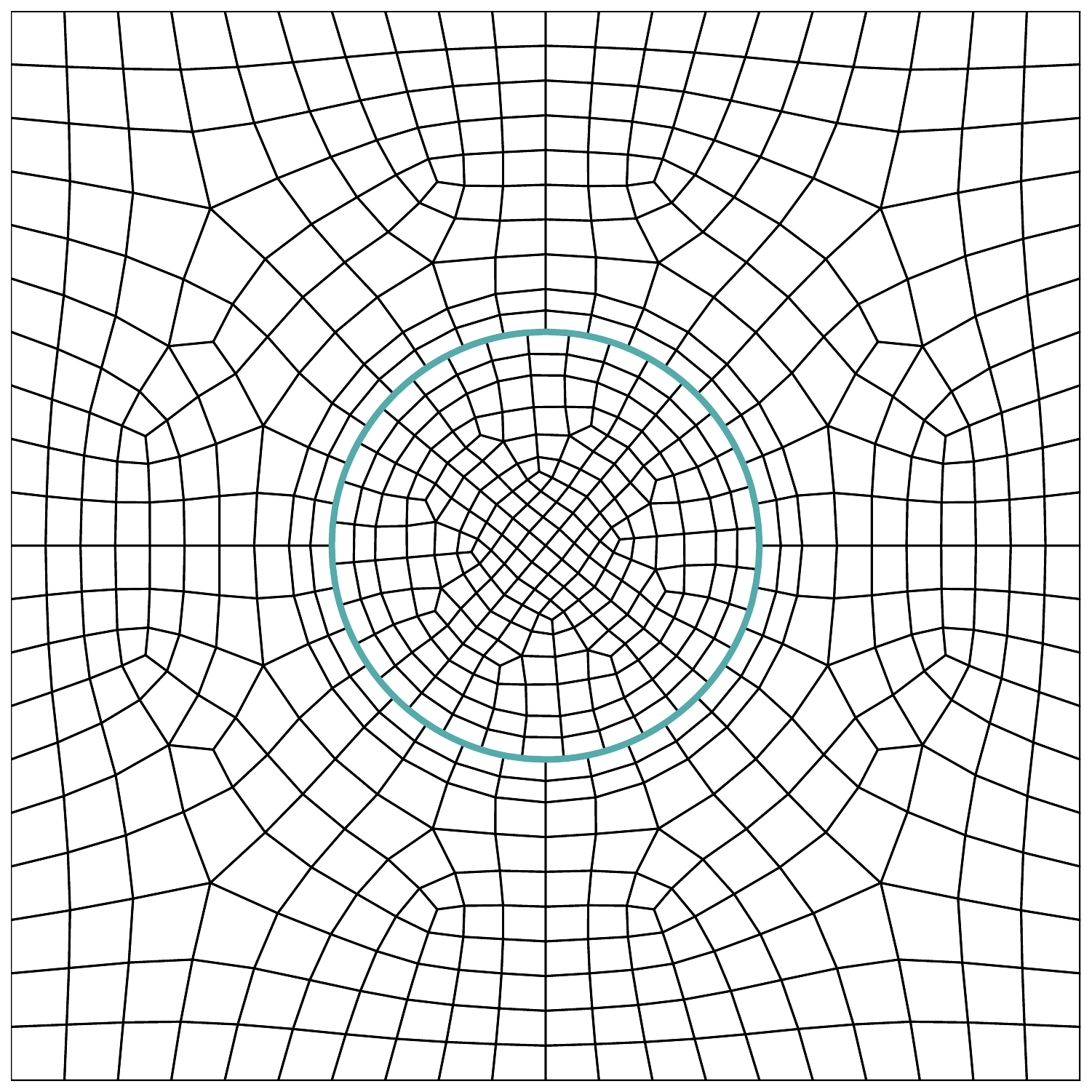}
  \caption{Mesh  with 700 cells at a time instant for the Euler vortex flow simulation (blue circle indicates sliding-mesh interface).}
  \label{fig:euler_grids}
\end{figure}

Figure (\ref{fig:euler_grids}) shows a computational mesh with 700 cells. The mesh has been decomposed into two parts: a rotating inner part with a radius of 2; a fixed outer part which takes the rest of the computational domain. Three meshes with 180, 700 and 2731 cells have been used for accuracy tests. For all three cases, the inner part is set to rotate at an angular speed of $\omega = 1.0$.

\begin{figure}[!htb]
  \vspace{0.7cm}
  \centering
  \includegraphics[scale=1.0]{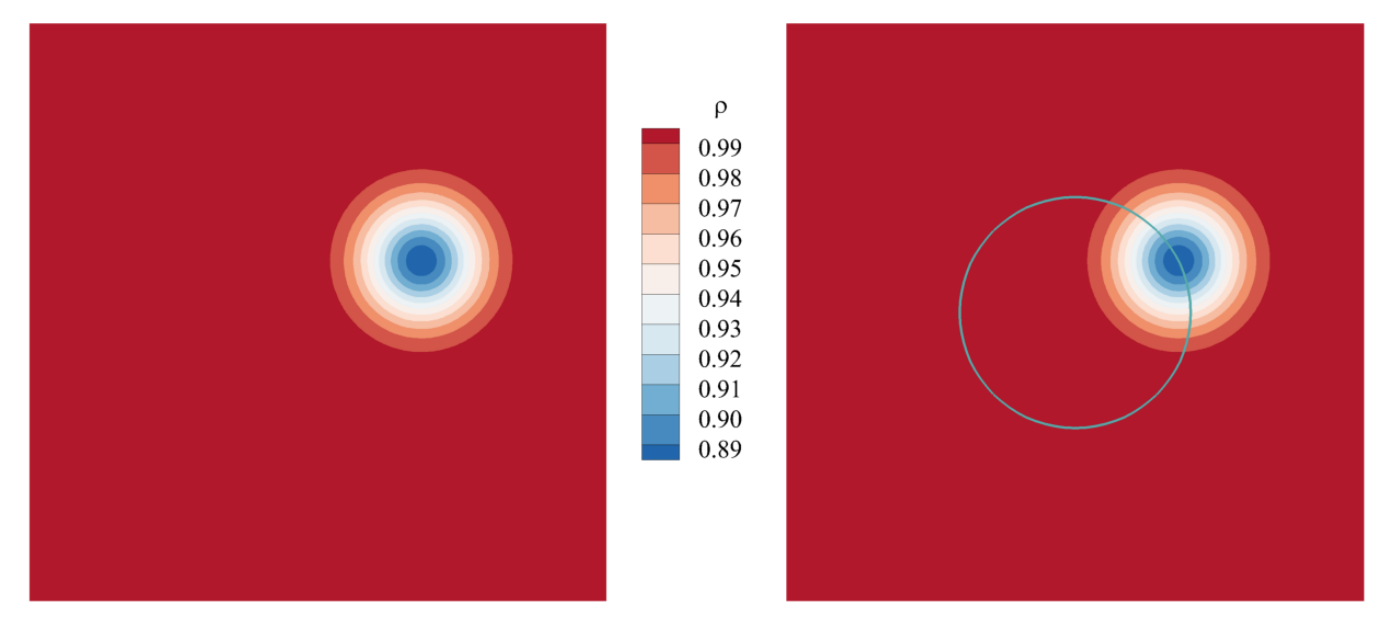} \\
  \caption{Contours of density at the time instant $t=2$ for the Euler vortex flow. Left, exact solution; right, numerical solution from fourth-order scheme (blue circle indicates location of sliding interface).}
  \label{fig:euler_rho}
\end{figure}

Figure (\ref{fig:euler_rho}) compares the density contours obtained with the fourth-oder SSD method on the finest mesh with the exact solution at $t=2$. As we can see, the solver resolves the vortex very well, and we see almost no visible difference between the exact solution and the numerical one.


Furthermore, Tables \ref{tab:euler3rd} and \ref{tab:euler4th} give the spatial accuracy of the scheme, where the $L1$ and $L2$ errors are computed from density at $t=2$ when vortex center is traveled right onto the sliding-mesh interface. From the two tables we see that the SSD method gives very reasonable order of accuracy.
\begin{table}[!ht]
  \vspace{0.7cm}
  \begin{center}
  \begin{tabularx}{0.7\textwidth}{XXXXl}
  \hline
  cells & L1 error &  order & L2 error & order \\
  \hline
  180   & 3.16E-4 &  -   & 8.13E-4 &  -    \\
  700   & 4.58E-5 & 2.85 & 1.13E-4 & 2.90  \\
  2731  & 7.00E-6 & 2.80 & 1.73E-5 & 2.83  \\
   \hline
  \end{tabularx}
  \end{center}
  \caption{Errors and orders of accuracy of the third-order scheme for the Euler vortex flow simulation.}
  \label{tab:euler3rd}
\end{table}
\begin{table}[!ht]
  \vspace{0.7cm}
  \begin{center}
  \begin{tabularx}{0.7\textwidth}{XXXXl}
  \hline
  cells & L1 error &  order & L2 error & order \\
  \hline
  180   & 5.43E-5 &  -   & 1.26E-4 &  -    \\
  700   & 3.27E-6 & 4.14 & 8.02E-6 & 4.06  \\
  2731  & 2.26E-7 & 4.03 & 5.50E-7 & 4.00  \\
  \hline
  \end{tabularx}
  \end{center}
  \caption{Errors and orders of accuracy of the fourth-order scheme for the Euler vortex flow simulation.}
  \label{tab:euler4th}
\end{table}

To see how efficient the SSD method is, we compare the total computational time and the communication time on the sliding-mesh interface in Tables \ref{tab:euler_eff_3rd} and \ref{tab:euler_eff_4th} for the third- and fourth-order schemes, respectively. Times in both tables are collected for 100 computational steps and do not include any post-processing time. It is seen that for all test cases, communication on the sliding-mesh interface takes only a few percent of the total computational time, which clearly shows that the SSD method is efficient. It is interesting that the relative communication time (represented by the percentage) decreases as either the number of cells or the order of schemes increases. This is due to the fact that cells are one dimension higher than faces: when we perform a mesh refinement, the total number of faces in the domain grows faster than on the sliding-mesh interface; when we increase the scheme order, the total number of degrees of freedom in the domain also grows faster than on the sliding-mesh interface.
\begin{table}[!ht]
    \vspace{0.7cm}
    \begin{center}
        \begin{tabularx}{0.7\textwidth}{XXXl}
            \hline
            cells & total time &  comm. time & percentage \\
            \hline
             180  & 0.254944 & 0.017657 & 6.92\% \\
             700  & 1.013039 & 0.041758 & 4.12\% \\
            2731  & 4.607205 & 0.097717 & 2.12\% \\
            \hline
        \end{tabularx}
    \end{center}
    \caption{Total computation time and interface communication time (both in seconds) for 100 computational steps using a third-order scheme for the Euler vortex flow simulation.}
    \label{tab:euler_eff_3rd}
\end{table}
\begin{table}[!ht]
    \vspace{0.7cm}
    \begin{center}
        \begin{tabularx}{0.7\textwidth}{XXXl}
            \hline
            cells & total time &  comm. time & percentage \\
            \hline
             180  & 0.420763 & 0.023282 & 5.53\% \\
             700  & 1.763656 & 0.058833 & 3.34\% \\
            2731  & 7.511551 & 0.125820 & 1.68\% \\
            \hline
        \end{tabularx}
    \end{center}
    \caption{Total computation time and interface communication time (both in seconds) for 100 computational steps using a fourth-order scheme for the Euler vortex flow simulation.}
    \label{tab:euler_eff_4th}
\end{table}

\subsection{Taylor-Couette flow}
\label{sec:couette}

\begin{figure}[!ht]
  \vspace{0.5cm}
  \centering
  \includegraphics[width=0.40 \textwidth, angle=0]{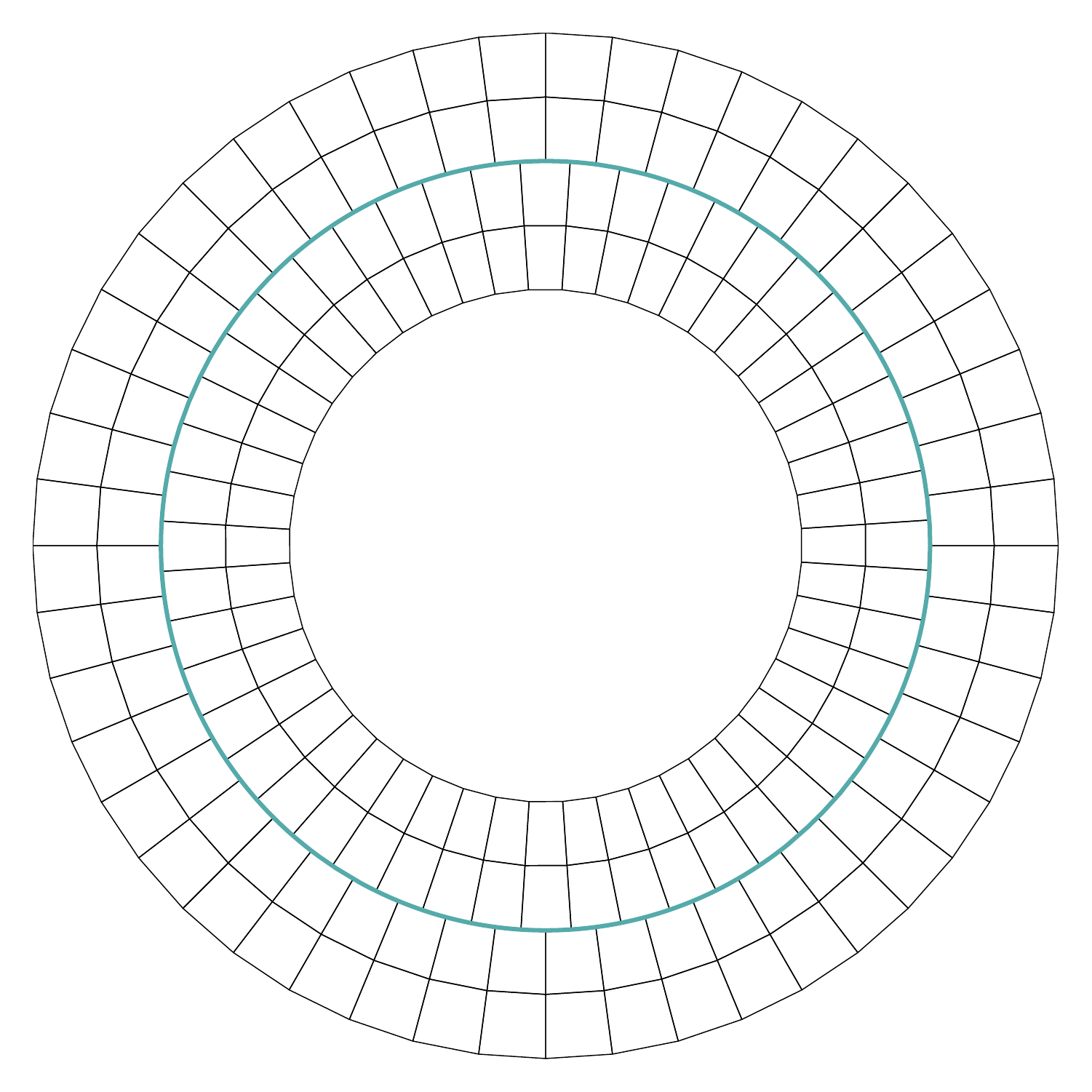}
  \caption{Mesh with 192 cells at a time instant for the Taylor-Couette flow simulation (blue circle indicates sliding-mesh interface).}
  \label{fig:couette_grids}
\end{figure}

To test the order accuracy on viscous flow, we use the laminar Taylor-Couette flow as the test case. Previous researchers such as \citet{liang-2009-cas}, \citet{michalak-09} used similar flows to test the accuracy of their solvers. In the present test, the inner cylinder has a radius of $r_i=1$, the outer cylinder has radius of $r_o=2$. Both boundaries are set to be isothermal walls. The domain has been divided into two parts at $r=1.5$. The inner part rotates at an angular speed of $\omega_i=1$ while the outer part stays stationary. The Reynolds number based the inner cylinder radius and speed is $Re=10$. The Mach number on the inner wall is set to $M_i=0.1$. Three meshes with 192, 768 and 3072 cells are used for the tests. Figure \ref{fig:couette_grids} shows the mesh with 192 cells.

\begin{figure}[!htb]
    \vspace{0.7cm}
    \centering
    \includegraphics[scale=1.0]{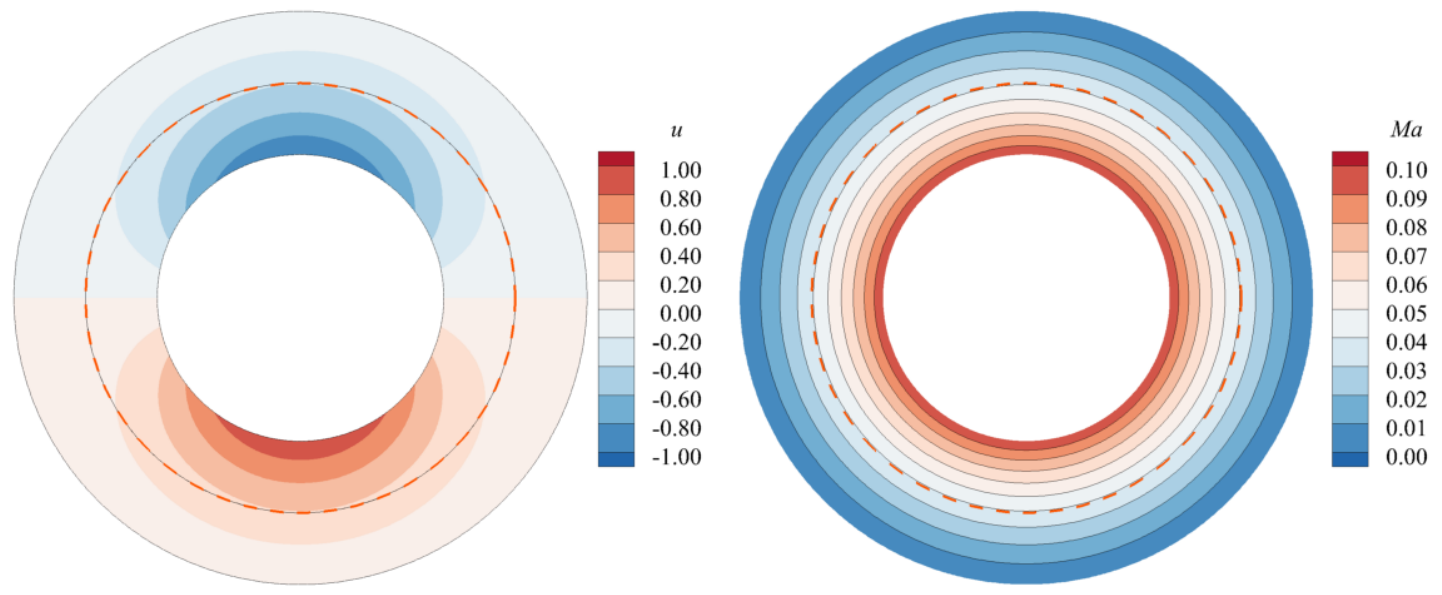} \\
    \caption{Contours of $u$ velocity component and the Mach number for the Taylor-Couette flow (dashed circle indicates location of sliding-mesh interface).}
    \label{fig:couette_contour}
\end{figure}

Figure \ref{fig:couette_contour} shows the steady state contours of the $u$ velocity component and the Mach number obtained with a fourth-oder scheme on the finest mesh. We see that the Mach contours at the  steady state are a series of concentric circles, and the $u$ velocity is highly symmetric. These results are consistent with our expectations.

The exact solution for the circumferential velocity has the following relation to radius $r$,
\begin{equation}
v_{\theta} = \omega_i r_i \frac{r_o/r - r/r_o}{r_o/r_i - r_i/r_o},
\end{equation}
The $x$ component of this velocity (i.e., $u$) is used to compute the $L1$ and $L2$ error norms. From Table \ref{tab:couette3rd} and Table \ref{tab:couette4th} we see that the SSD method preserves the high-order accuracy for viscous flow as well.

\begin{table}[!ht]
  \vspace{0.7cm}
  \begin{center}
  \begin{tabularx}{0.7\textwidth}{XXXXl}
  \hline
  cells & L1 error &  order & L2 error & order \\
  \hline
  192  & 5.90E-5 &  -   & 8.71E-5 &  -    \\
  768  & 6.95E-6 & 3.09 & 9.82E-6 & 3.15  \\
  3072 & 8.43E-7 & 3.07 & 1.09E-6 & 3.16  \\
  \hline
  \end{tabularx}
  \end{center}
  \caption{Errors and orders of accuracy of a third-order scheme for the Taylor-Couette flow simulation.}
  \label{tab:couette3rd}
\end{table}
\begin{table}[!ht]
  \vspace{0.7cm}
  \begin{center}
  \begin{tabularx}{0.7\textwidth}{XXXXl}
  \hline
  cells & L1 error &  order & L2 error & order \\
  \hline
  192  & 9.72E-6 &  -   & 1.52E-5 &  -    \\
  768  & 6.16E-7 & 3.98 & 1.01E-6 & 3.91  \\
  3072 & 4.22E-8 & 3.92 & 6.60E-8 & 3.92  \\
  \hline
  \end{tabularx}
  \end{center}
  \caption{Errors and orders of accuracy of a fourth-order scheme for the Taylor-Couette flow simulation.}
  \label{tab:couette4th}
\end{table}

The total computational time and the sliding-mesh interface communication time are shown in Table \ref{tab:couette_eff_3rd} and Table \ref{tab:couette_eff_4th} for third- and fourth-order schemes, respectively. Again, data in both tables is collected for 100 computational steps and do not include any post-processing time. It is seen that for viscous flow simulations the SSD method remains very efficient.
\begin{table}[!ht]
    \vspace{0.7cm}
    \begin{center}
        \begin{tabularx}{0.7\textwidth}{XXXl}
            \hline
            cells & total time &  comm. time & percentage \\
            \hline
            192   & 1.267700  & 0.105698 & 8.34\% \\
            768   & 4.394639  & 0.234498 & 5.34\% \\
            3072  & 18.138023 & 0.561518 & 3.10\% \\
            \hline
        \end{tabularx}
    \end{center}
    \caption{Total computation time and interface communication time (both in seconds) for 100 computational steps using a third-order scheme for the Taylor-Couette flow simulation.}
    \label{tab:couette_eff_3rd}
\end{table}
\begin{table}[!ht]
    \vspace{0.7cm}
    \begin{center}
        \begin{tabularx}{0.7\textwidth}{XXXl}
            \hline
            cells & total time &  comm. time & percentage \\
            \hline
            192   & 2.075053  & 0.139513 & 6.92\% \\
            768   & 7.628162  & 0.322616 & 4.12\% \\
            3072  & 30.852314 & 0.729768 & 2.12\% \\
            \hline
        \end{tabularx}
    \end{center}
    \caption{Total computation time and interface communication time (both in seconds) for 100 computational steps using a fourth-order scheme for the Taylor-Couette flow simulation.}
    \label{tab:couette_eff_4th}
\end{table}

\subsection{Flow over a rotating elliptic cylinder}
To further verify the approach for flow involving complex geometries, we simulate a laminar flow over a two-dimensional elliptic cylinder in this section. \citet{Maruoka-03} and \citet{XZhang-08} studied incompressible flow over a rotating elliptic cylinder using the finite element and finite volume methods, respectively. Both studies use Chimera grids for communication between the foreground rotating mesh and the background stationary mesh. The freestream Mach number in the present test is set to $M_\infty=0.05$ to mitigate the compressibility effects in order to compare with their incompressible flow results.

\begin{figure}[!ht]
    \centering
    \includegraphics[scale=0.9]{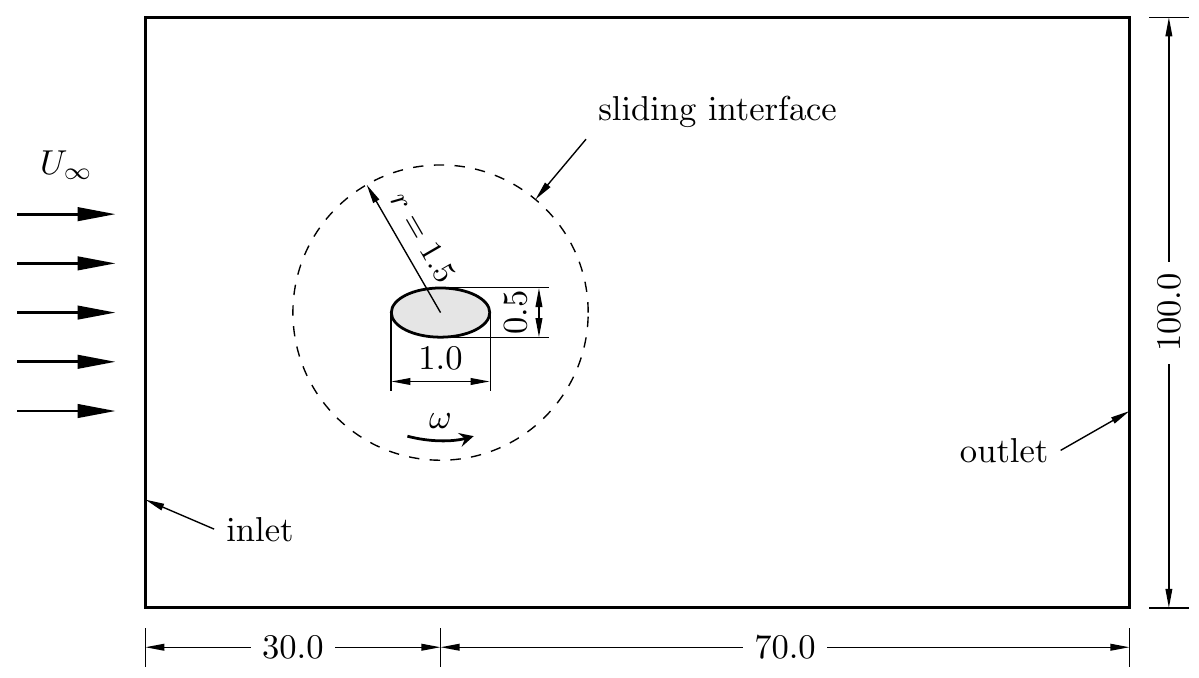} \\
    \caption{Schematic of the computational domain for flow over a rotating elliptic cylinder (not to scale).}
    \label{fig:cyl_domain}
\end{figure}

The elliptic cylinder has a major and a minor axis lengths of 1.0 and 0.5, respectively. Initially, the major axis is parallel to the freestream. The cylinder rotates counterclockwisely at an angular speed of $\omega=0.5\pi$. The Reynolds number based on the freestream velocity and the major axis length is $Re=200$. Figure \ref{fig:cyl_domain} shows a schematic of the computational domain. Slip boundary conditions are applied on the top and bottom boundaries. Dirichlet boundary condition is used at the inlet, and fixed pressure boundary condition is applied at the outlet. Finally, no-slip isothermal wall boundary condition is used on the cylinder surface.
\vspace{0.5cm}

\begin{figure}[!htb]
    \vspace{0.5cm}
    \centering
    \includegraphics[height=0.28 \textwidth, angle=0]{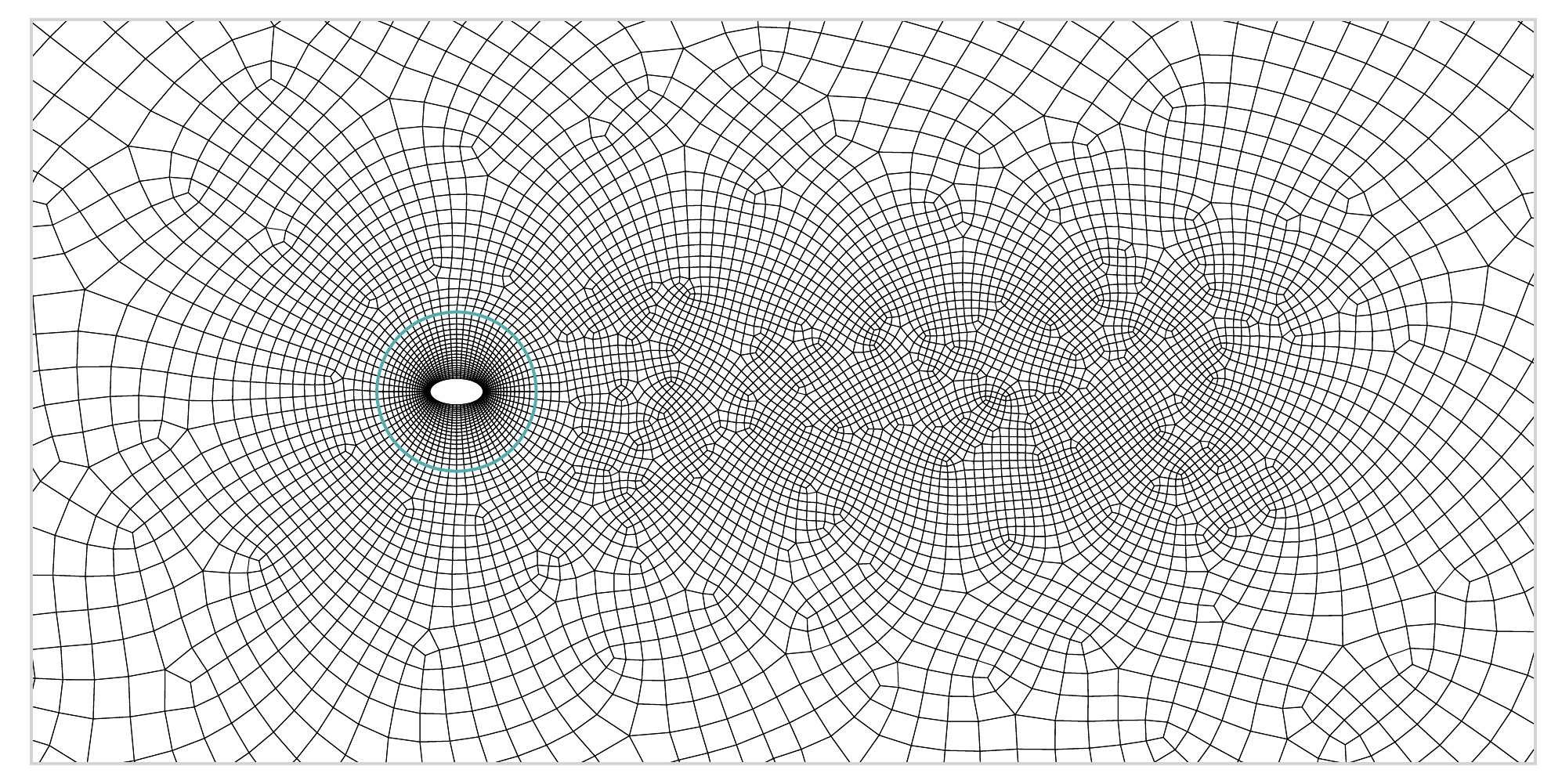}
    \includegraphics[height=0.28 \textwidth, angle=0]{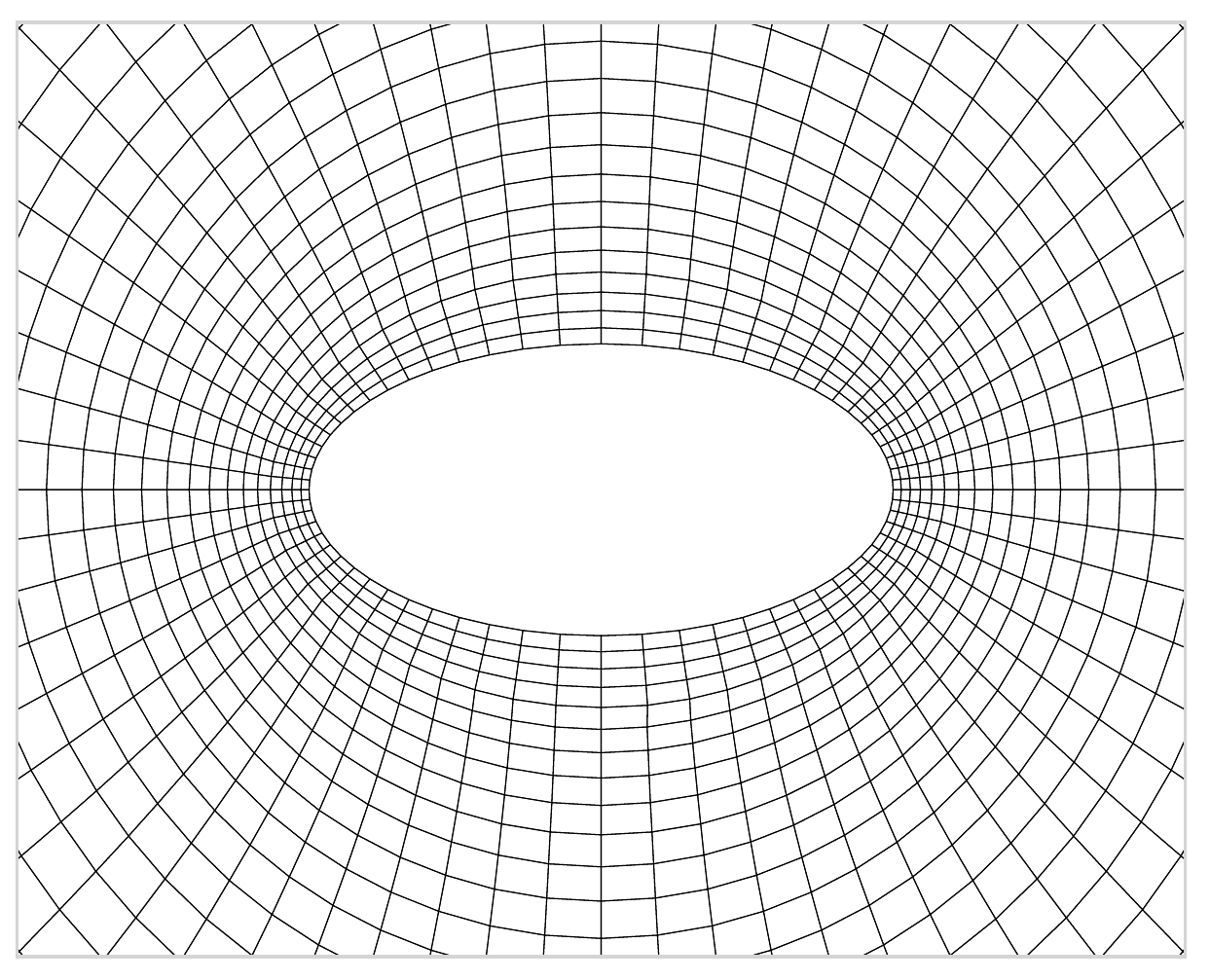}\\
    \caption{Two local views of the mesh for simulation of flow around a rotating elliptic cylinder (blue circle indicates sliding-mesh interface).}
    \label{fig:cyl_msh}
\end{figure}

The inner rotating domain has a radius of 1.5 and is meshed with 1280 cells. The rest of the domain is stationary and has 7391 cells. Mesh refinement are performed around the leading and the trailing edges, as well as in the wake region. Figure \ref{fig:cyl_msh} shows two local views of the mesh. The nondimensional time step size $\Delta t  U_{\infty}/L$ for the simulation is set to $1.0\times 10^{-4}$, where L is the major axis length and $U_{\infty}$ is the freestream velocity.

\begin{figure}[!htb]
    \vspace{0.5cm}
    \centering
    \includegraphics[scale=1.0]{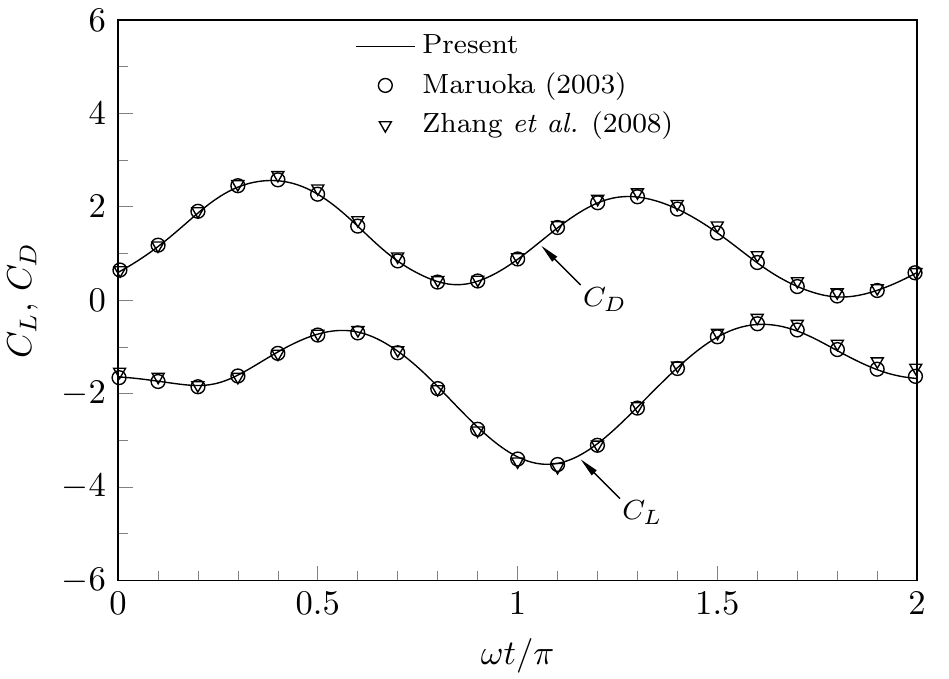}\\
    \caption{Lift and drag coefficients for flow over an elliptic cylinder.}
    \label{fig:cyl_clcd}
\end{figure}

Both third- and fourth-order schemes were tested for this flow and no visible difference was observed between the solutions. For this reason, we only present results from the fourth-order scheme. As was noticed by \citet{Maruoka-03} and \citet{XZhang-08}, the fully developed flow takes a periodic pattern as the cylinder rotates. The lift and the drag coefficients in one rotating period are shown in Figure \ref{fig:cyl_clcd}. It is seen that the present results agree very well with the previously published results.

\begin{figure}[!htb]
    \vspace{0.5cm}
    \centering
    \includegraphics[width=0.92 \textwidth]{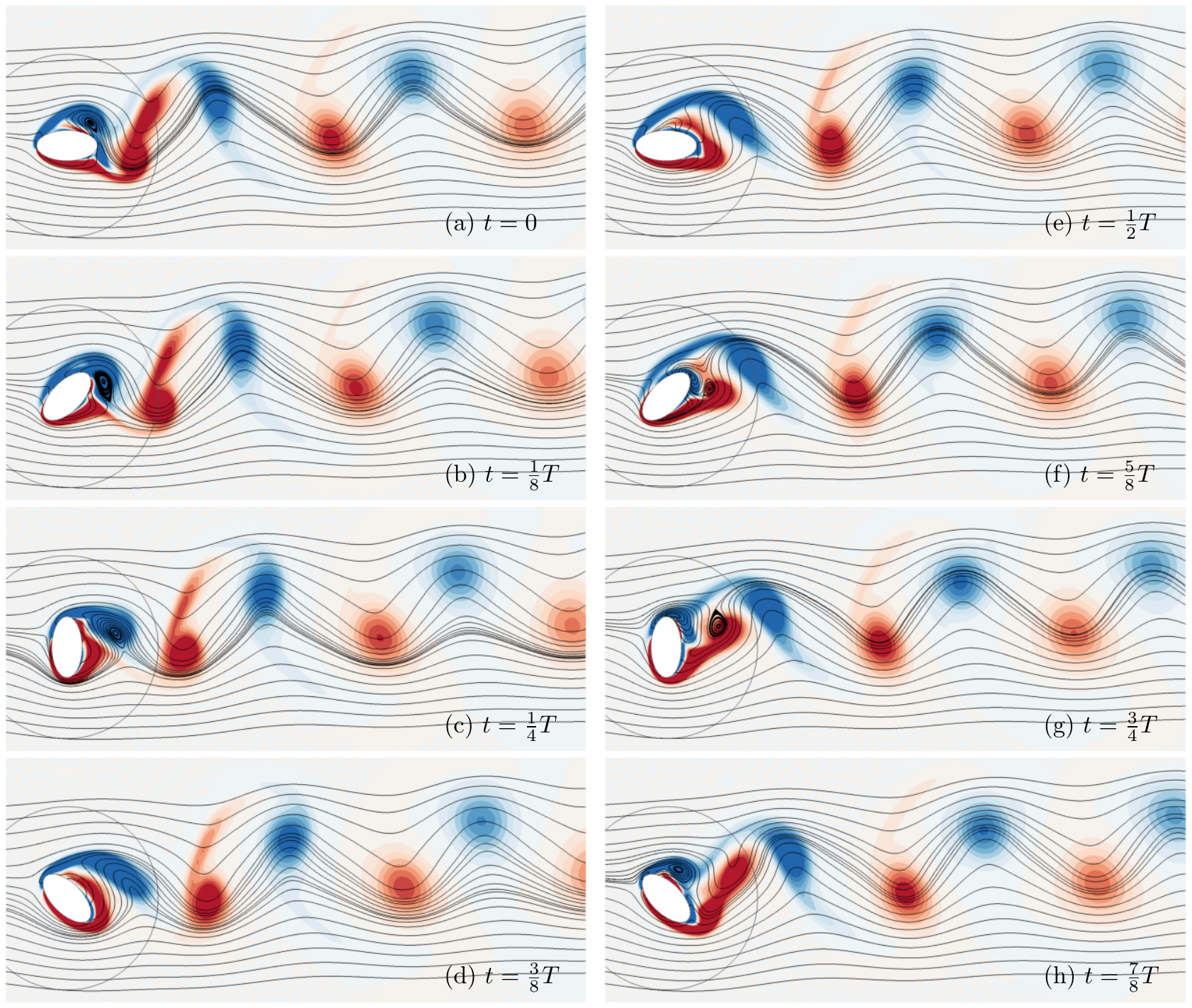}\\
    \caption{Streamlines and vorticity contours (blue color means negative value, red color means positive) for flow over a rotating elliptic cylinder.}
    \label{fig:cyl_flow}
\end{figure}

Figure \ref{fig:cyl_flow} shows the streamlines superimposed on vorticity contours at a series of time instants in one rotating period. A clockwise and a counterclockwise vortices appear alternatively around the two ends of the cylinder. From Figure \ref{fig:cyl_flow} (h) and Figure \ref{fig:cyl_flow} (a), we see that a clockwise vortex is formed at the cylinder leading edge as the cylinder rotates, and this vortex is then shed off from the leading edge and travels downstream to hit the trailing edge. From Figure \ref{fig:cyl_flow} (b) to Figure \ref{fig:cyl_flow} (d), the same clockwise vortex is again shed off from the trailing edge and then convected towards downstream. From Figure \ref{fig:cyl_flow} (e) to Figure \ref{fig:cyl_flow} (g), a counterclockwise vortex slowly emerges around the other end of the cylinder and is then convected downstream without reattaching to the cylinder. This process repeats as the cylinder rotates, and a vortex street forms downstream of the cylinder.

The efficiency of the SSD method is shown in Table \ref{tab:cyl_eff} for this case. The results in the table confirm our previous conclusion, i.e., the relative communication time generally decreases as the number of cells or the order of scheme increases. In fact, the interface communication time in the table is almost negligible comparing to the total computational time.

\begin{table}[!ht]
    \vspace{0.4cm}
    \begin{center}
        \begin{tabularx}{0.7\textwidth}{XXXl}
            \hline
            order & total time &  comm. time & percentage \\
            \hline
            3   & 38.924785  & 0.143083 & 0.37\% \\
            4   & 70.655175  & 0.185890 & 0.26\% \\
            \hline
        \end{tabularx}
    \end{center}
    \caption{Total computation time and interface communication time (both in seconds) for 100 computational steps for simulation of flow over a rotating elliptic cylinder.}
    \label{tab:cyl_eff}
\end{table}

\subsection{Flow inside a 2D stirred tank}
\begin{figure}[!ht]
    \centering
    \includegraphics[width=0.5 \textwidth, angle=0]{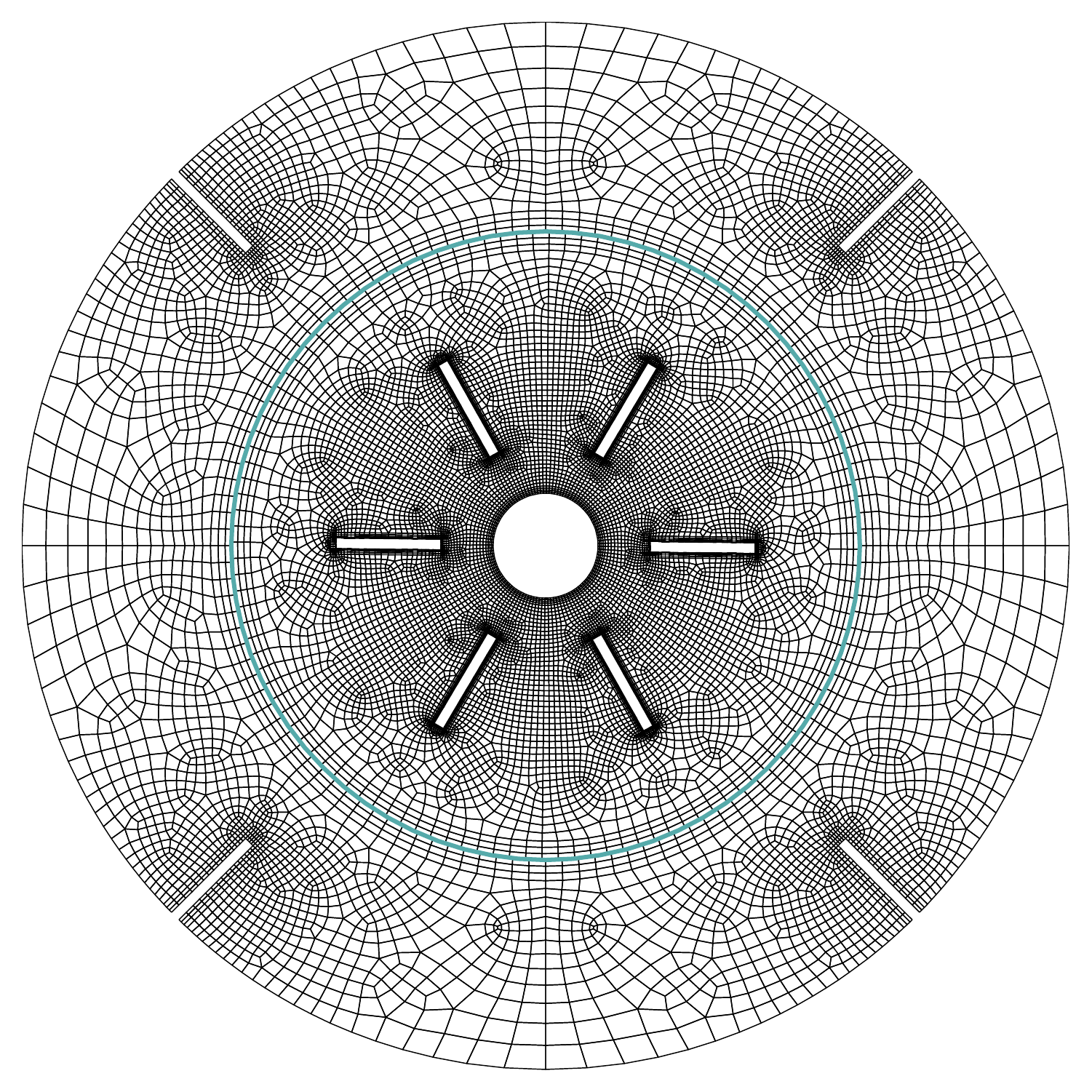}
    \caption{Mesh for simulation of flow inside a $2D$ stirred tank (blue circle indicates sliding-mesh interface).}
    \label{fig:stirredtank_grids}
\end{figure}

In this last test case, we apply the SSD method to simulate laminar flow inside a $2D$ stirred tank. Figure \ref{fig:stirredtank_grids} shows the unstructured quadrilateral mesh used for this case. As we can see, the tank has several components: an inner circular wall with a radius of 0.5; an outer circular wall with a radius of 5; six uniformly distributed agitating blades, each extends from $r=1$ to $r=2$ and has a thickness of 0.1; four baffles installed on the outer wall, each of them has a height of 1 and a thickness of 0.1. The computational domain is split into an inner rotating part and an outer fixed part, resulting in a sliding-mesh interface at $r=3$. Mesh has been refined around the blades, baffles, and the wall boundaries. The resulted mesh has $14,990$ cells in total.

In this simulation, the inner circular wall and the blades rotate at an angular speed of $\omega=1$. The Reynolds number base on the inner wall diameter and the angular speed is $Re=100$. Flow on the inner wall surface has a Mach number of $M_i=0.1$. No slip isothermal wall boundary conditions are used on the inner wall surface, the baffles, and the outer walls. Adiabatic wall boundary conditions are used on the six blades. A time step size of $\Delta t=1.0\times 10^{-4}$ is used.

\begin{figure}[!ht]
  \vspace{0.5cm}
  \centering
  \includegraphics[scale=1.0]{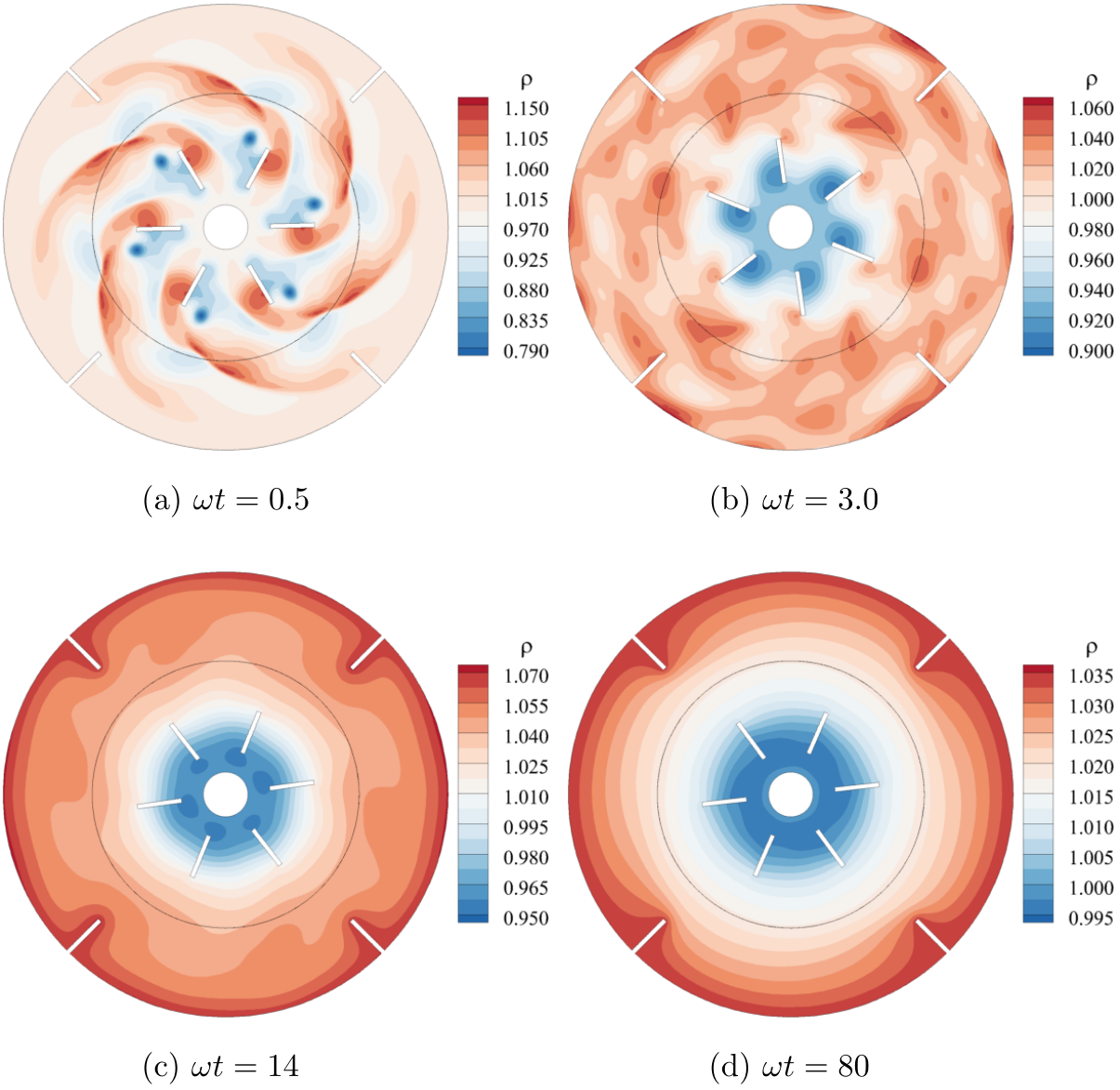} 
  \caption{Density contours of flow inside a $2D$ stirred tank at four distinct states.}
  \label{fig:stirredtank_rho}
\end{figure}

Initially, the flow field is set to be uniform and stationary. Figure \ref{fig:stirredtank_rho} shows the flow fields at four different states by visualizing the density. It is seen in Figure \ref{fig:stirredtank_rho}(a) that at the initial transient state, the flow behaves like flow around a propeller: fluid is ``squeezed" and ``pushed" away from the blades, resulting in a concentric flow pattern, and large fluctuations are observed in the flow field. Quickly after the flow reaches the baffles and the outer wall, the flow field becomes very chaotic: vortical structures generated by flow passing the baffles, bouncing pressure waves, blades induced vortices, unsteady boundary layers, etc., as can be seen in Figure \ref{fig:stirredtank_rho}(b). As the blades continue rotating for a longer time, the chaotic flow structures are slowly being dissipated, and organized large flow structures emerge as shown in Figure \ref{fig:stirredtank_rho}(c). Finally the flow field reaches a quasi steady state in the rotating reference frame, carrying little variation with time, as shown in Figure \ref{fig:stirredtank_rho}(d). Interestingly, the flow is becoming more and more uniform at the later stage of stirring. Variations on the density field get smaller and smaller when comparing Figure \ref{fig:stirredtank_rho} (a)-(d) along time.

As for all previous cases, we monitored the efficiency of the SSD method for this case. The results are shown in Table \ref{tab:tank_eff}. Again, as can be seen from the table, the SSD method is very efficient with negligible communication time.
\begin{table}[!ht]
    \vspace{0.4cm}
    \begin{center}
        \begin{tabularx}{0.7\textwidth}{XXXl}
            \hline
            order & total time &  comm. time & percentage \\
            \hline
            3   & 85.547403   & 0.328987 & 0.38\% \\
            4   & 148.780491  & 0.437902 & 0.29\% \\
            \hline
        \end{tabularx}
    \end{center}
    \caption{Total computation time and interface communication time (both in seconds) for 100 computational steps for simulation of flow inside a $2D$ stirred tank.}
    \label{tab:tank_eff}
\end{table}

\section{Conclusions}
\label{sec:conclusion}
In this paper, a novel, simple, efficient, and high-order accurate sliding-mesh interface method is reported for subsonic compressible flows. The sliding-mesh spectral difference (SSD) method has been successfully developed and tested for several inviscid and viscous flow problems. The SSD method retains the high-order accuracy of the SD method. It is also shown that the SSD method is very efficient as it introduces negligible extra computational cost to the SD method for realistic flow simulations. In this paper, we demonstrated the approach on uniformly meshed sliding-mesh interfaces where each cell face has two mortars. The sliding-mesh interface method can be extended to handle more than two mortars for each face. The SSD method is also very suitable for parallel computing. Since there is no overlapping in grids and the sliding-mesh interface method introduces negligible computational time, thus each domain can be decomposed and distributed to the processors to achieve load balancing. Finally, this high-order curved sliding-mesh interface method can also be extended to other discontinuous high-order methods for compressible flows.

\section{Acknowledgment}
The authors would like to express our acknowledgments for an Office of Naval Research (ONR) support with award number N000141210500 managed by Dr. Ki-Han Kim. Chunlei Liang would like to also acknowledge the support by an ONR Young Investigator Program award administrated by Dr. Ki-Han Kim.




\newpage
\bibliographystyle{elsarticle-harv}
\bibliography{reference}

\newpage
\listoftables

\newpage
\listoffigures

\end{document}